\pdfoutput=1
\RequirePackage{ifpdf}
\ifpdf 
\documentclass[pdftex]{sigma}
\else
\documentclass{sigma}
\fi

\DeclareMathOperator{\Tr}{Tr}

\begin{document}

\newcommand{\res}[1]{\underset{#1}{\mathrm{Res}}}

\newcommand{\limr}[2]{\lim\limits_{#1 \to #2^{+}} (#1 - #2)}

\newcommand{\arXivNumber}{1405.5396}

\allowdisplaybreaks

\renewcommand{\PaperNumber}{097}

\FirstPageHeading

\ShortArticleName{Quantum Dimension and Quantum Projective Spaces}

\ArticleName{Quantum Dimension and Quantum Projective Spaces}

\Author{Marco MATASSA}

\AuthorNameForHeading{M.~Matassa}

\Address{SISSA, Via Bonomea 265, I-34136 Trieste, Italy}
\Email{\href{mailto:marco.matassa@gmail.com}{marco.matassa@gmail.com}}

\ArticleDates{Received July 29, 2014, in f\/inal form September 21, 2014; Published online September 25, 2014}

\Abstract{We show that the family of spectral triples for quantum projective spaces introduced by D'Andrea and
D\c{a}browski, which have spectral dimension equal to zero, can be reconsidered as modular spectral triples by taking
into account the action of the element $K_{2\rho}$ or its inverse.
The spectral dimension computed in this sense coincides with the dimension of the classical projective spaces.
The connection with the well known notion of quantum dimension of quantum group theory is pointed out.}

\Keywords{quantum projective spaces; quantum dimension; modular spectral triples}

\Classification{58J42; 58B32; 46L87}

\section{Introduction}

Quantum homogeneous spaces provide an excellent testing ground to study how quantum groups f\/it into the framework of
non-commutative geometry developed by Connes~\cite{con-book}.
An important result in this respect is given in~\cite{qflag}, where a~Dirac operator~$D$ is def\/ined on quantized
irreducible generalized f\/lag manifolds, which yields a~Hilbert space realization of the covariant f\/irst-order
dif\/ferential calculus constructed in~\cite{flag-calc}.
This in particular means that the commutator of~$D$ with an element of the coordinate algebra is a~bounded operator,
which is one of the def\/ining properties of a~spectral triple.
On the other hand the compactness of the resolvent of~$D$, which is another essential requirement of the theory, has not
been proven yet, even though it is expected to hold.
In particular it can be checked for the simplest case of this construction, the Podle\'{s} sphere, in which case~$D$
coincides with the Dirac operator introduced in~\cite{DS-pod}, which has compact resolvent.
In this respect we should also mention the work~\cite{krtu13}, where Dolbeault--Dirac operators are constructed for
Hermitian symmetric spaces, which can be seen as a~f\/irst step towards checking the compactness condition.

Among the class of~$q$-deformed irreducible f\/lag manifolds are the quantum projective spa\-ces~$\mathbb{C}{\rm P}_{q}^{\ell}$,
the simplest example of which is again the Podle\'{s} sphere, which is obtained for~$\ell=1$.
The case of~$\mathbb{C}{\rm P}_{q}^{2}$ has been studied in~\cite{projectiveplane} and then generalized in~\cite{projective}
to $\mathbb{C}{\rm P}_{q}^{\ell}$ with $\ell \geq 2$.
The starting point is the introduction of the~$q$-analogue of the module of antiholomorphic dif\/ferential $k$-forms~$\Omega^{k}$.
More generally the modules $\Omega^{k}_{N}$ are considered, with $N \in \mathbb{Z}$, corresponding essentially to
$\Omega^{k} = \Omega^{k}_{0}$ twisted by certain line bundles.
The Hilbert space completion of these is denoted by $H_{N}$.
For each of these an unbounded self-adjoint operator $D_{N}$ is introduced, which has bounded commutators with the
coordinate algebra $\mathcal{A}(\mathbb{C}{\rm P}_{q}^{\ell})$.
The main result is that $(\mathcal{A}(\mathbb{C}{\rm P}_{q}^{\ell}), H_{N}, D_{N})$ is a~family of equivariant spectral
triples.

It turns out that these spectral triples are $0^{+}$-summable, in the sense that the operator
$(D_{N}^{2}+1)^{-\epsilon/2}$ is trace-class for every $\epsilon > 0$.
The detailed computation of the spectrum clearly reveals why this is the case: the eigenvalues of this operator grow
like a~$q$-number, so exponentially, while their multiplicities grow like a~polynomial.
We recall that in the classical case it is the balance between the growth of the eigenvalues and their multiplicities
that allows to recover the dimension of the manifold in consideration.
In this case the eigenvalues grow much faster than their multiplicities, which explains the $0^{+}$-summability.

In this paper we consider a~simple modif\/ication of the above construction, which f\/its into the framework of modular
spectral triples.
The idea is to consider the action of the element $K_{2\rho}$, which implements the modular group of the Haar state of
$\mathcal{A}(\mathbb{C}{\rm P}_{q}^{\ell})$.
In particular we compute the spectral dimension associated to~$D$ with respect to the weight $\Tr(K_{2\rho}
\cdot)$, with the result that it coincides with the classical dimension.
This computation is linked with an important concept in the theory of quantum groups, which is the notion of quantum
dimension.
We also point out that, as a~consequence of a~property of the quantum dimension, the same result for the spectral
dimension is obtained by considering $K_{2\rho}^{-1}$.
This in turn is connected with results from twisted Hochschild (co)homology, as we will discuss in the last part.

One of the motivations for this paper comes from the notion of integration which is def\/ined in the context of spectral
triples.
Let us brief\/ly review how this works in the classical case, by considering the canonical spectral triple
$(C^{\infty}(M), L^{2}(M,S), D)$ associated to a~compact spin manifold~$M$.
First of all, the dimension of such a~manifold can be extracted from the spectrum of~$D$.
Indeed the operator $(D^{2}+1)^{-s/2}$, with $s \in \mathbb{R}$, turns out to be trace-class for all $s > n$, where~$n$
is the dimension of~$M$.
The operator~$D$ also allows to def\/ine a~notion of integration.
One possible formulation is via the linear functional $\psi: C^{\infty}(M) \to \mathbb{C}$ given~by
\begin{gather*}
\psi(f) = \limr{s}{n} \Tr\big(f \big(D^{2}+1\big)^{-s/2}\big),
\end{gather*}
where~$f$ is considered to be acting via left multiplication.
This functional turns out to coincide with the usual integral of~$f$, up to a~constant.
Therefore this procedure allows to recover the dimension and the integration of functions for the manifold in
consideration.

As we mentioned above, in the case of quantum projective spaces the analogue of this procedure gives a~spectral
dimension equal to zero.
But, more importantly, it does not allow to recover the natural notion of integration that is available on these spaces,
which is given by the faithful Haar state.
The reason for this failure is quite clear: the analogue of the functional~$\psi$ is by construction a~trace, as it
follows from the def\/ining properties of a~spectral triple, while on the other hand the Haar state is not.
A~way out of this problem is to replace the trace in the def\/inition of~$\psi$ by a~weight, which immediately brings us
into the realm of modular spectral triples.
Similar ideas were discussed in~\cite{mmintegration} for the case of ${\rm SU}_{q}(2)$.

The plan of the paper is as follows.
In Section~\ref{sec:projective} we recall some basic notions on quantum projective spaces and their family of spectral
triples.
In Section~\ref{sec:modular} we recall the notion of modular spectral triple and discuss its role in the context of this
paper.
In Section~\ref{sec:quantum} we show how the computation of the spectral dimension is connected to the notion of quantum
dimension, and compute it in the cases of interest to us.
Finally in Section~\ref{sec:spectral} we prove that the spectral dimension coincides with the dimension of the classical
projective spaces.
We also discuss the connection with Hochschild twisted (co)homology.

\section{Quantum projective spaces}
\label{sec:projective}

In this section we provide some background on quantum projective spaces, which we denote by $\mathbb{C}{\rm P}_{q}^{\ell}$ for
$\ell \in \mathbb{N}$ and $\ell \geq 2$.
These are~$q$-deformations of complex projective spaces of real dimension $2\ell$.
The case $\ell = 1$ of this construction coincides with the standard Podle\'{s} sphere and is well known in the
literature.
We take our def\/initions and notations from~\cite{projective}.

To def\/ine quantum projective spaces we f\/irst def\/ine the Hopf $*$-algebra ${\rm U}_{q}(\mathfrak{su}(\ell + 1))$, which is
a~deformation of the universal enveloping algebra ${\rm U}(\mathfrak{su}(\ell + 1))$, and its dual $\mathcal{A}({\rm SU}_{q}(\ell +
1))$, which can be considered as the algebra of representative functions on the quantum ${\rm SU}(\ell+1)$ group.
Our reference for this material is the book~\cite{KS}, but one must keep in mind that what we denote~by
${\rm U}_{q}(\mathfrak{g})$ corresponds to $\breve{\rm U}_{q}(\mathfrak{g})$ there.
The coordinate algebra $\mathcal{A}(\mathbb{C}{\rm P}_{q}^{\ell})$ of the quantum projective space $\mathbb{C}{\rm P}_{q}^{\ell}$
can then be def\/ined as the f\/ixed point subalgebra of $\mathcal{A}({\rm SU}_{q}(\ell + 1))$ for the action of a~suitable Hopf
subalgebra of ${\rm U}_{q}(\mathfrak{su}(\ell + 1))$.
We now review these notions.

For $0 < q < 1$ we denote by ${\rm U}_{q}(\mathfrak{su}(\ell + 1))$ the $*$-algebra generated by $K_{i} = K_{i}^{*}$,
$K_{i}^{-1}$, $E_{i}$ and $F_{i} = E_{i}^{*}$, with $i = 1, \dots, \ell$, and with relations
\begin{gather*}
[K_{i}, K_{j}] = 0,
\qquad
K_{i} E_{i} K_{i}^{-1} = q E_{i},
\\
K_{i} E_{j} K_{i}^{-1} = q^{-1/2} E_{j}
\qquad
\text{if}\quad |i-j| = 1,
\\
K_{i} E_{j} K_{i}^{-1} = E_{j}
\qquad
\text{if}\quad |i-j| > 1,
\\
[E_{i}, F_{j}] = \delta_{ij} \frac{K_{i}^{2} - K_{i}^{-2}}{q-q^{-1}},
\\
E_{i}^{2}E_{j} - \big(q+q^{-1}\big)E_{i}E_{j}E_{i} + E_{j}E_{i}^{2} = 0
\qquad
\text{if}\quad |i-j| = 1,
\\
[E_{i},E_{j}] = 0
\qquad
\text{if}\quad |i-j| > 1.
\end{gather*}
We call ${\rm U}_{q}(\mathfrak{su}(\ell))$ the Hopf $*$-subalgebra generated by the elements $K_{i} = K_{i}^{*}$,
$K_{i}^{-1}$, $E_{i}$ and $F_{i} = E_{i}^{*}$ with $i = 1, \dots, \ell - 1$.
Its commutant is the Hopf $*$-subalgebra ${\rm U}_{q}(\mathfrak{u}(1))$ generated by the element $K_{1} K^{2}_{2} \cdots
K^{\ell}_{\ell}$ and its inverse.
This is a~positive operator in all the representations we consider, so we can def\/ine its root of order $2/(\ell+1)$~by
\begin{gather*}
\hat{K} = \big(K_{1} K^{2}_{2} \cdots K^{\ell}_{\ell}\big)^{2/(\ell+1)}.
\end{gather*}

The following element will play a~central role in this paper:
\begin{gather*}
K_{2\rho} = \big(K_{1}^{\ell} K_{2}^{2(\ell-1)} \cdots K_{j}^{j(\ell-j+1)} \cdots K_{\ell}^{\ell}\big)^{2}.
\end{gather*}
Here the symbol~$\rho$ denotes the Weyl vector of the Lie algebra $\mathfrak{su}(\ell+1)$, see for example~\cite{KS} for
its role in~$q$-deformations of semisimple Lie algebras.
One important property of this element is that it implements the square of the antipode, in the sense that $S^{2}(h) =
K_{2\rho} h K_{2\rho}^{-1}$ for any $h \in {\rm U}_{q}(\mathfrak{su}(\ell + 1))$.
More importantly for us, it also implements the modular group of the Haar state of $\mathcal{A}({\rm SU}_{q}(\ell + 1))$, as
we will see in a~moment.

We are interested in representations in which $K_{j}$ is represented by a~positive operator.
Such irreducible f\/inite-dimensional $*$-representations of ${\rm U}_{q}(\mathfrak{su}(\ell + 1))$ are labeled by~$\ell$
non-negative integers.
Writing $n = (n_{1}, \dots, n_{\ell})$, we denote by $V_{n}$ the vector space carrying the representa\-tion~$\rho_{n}$
with highest weight~$n$.
This means that there exists a~vector~$v$ which is annihilated by all the $E_{j}$'s and satisf\/ies $\rho_{n}(K_{i})v =
q^{n_{i}/2}v$.

We now introduce the coordinate algebra $\mathcal{A}({\rm SU}_{q}(\ell + 1))$.
It is the Hopf $*$-algebra generated by the elements $u^{i}_{j}$, with $i,j = 1, \dots, \ell+1$, and with relations
\begin{gather*}
u^{i}_{k} u^{j}_{k} = q u^{j}_{k} u^{i}_{k},
\qquad
u^{k}_{i} u^{k}_{j} = q u^{k}_{j} u^{k}_{i},
\qquad
\text{for}\quad i < j,
\\
\big[u^{i}_{l},u^{j}_{k}\big] = 0,
\qquad
\big[u^{i}_{k},u^{j}_{l}\big] = \big(q-q^{-1}\big) u^{i}_{l}u^{j}_{k},
\qquad
\text{for}\quad i < j, \quad k < l.
\end{gather*}
and with the determinant relation
\begin{gather*}
\sum\limits_{p \in S_{\ell+1}}(-q)^{\|p\|} u^{1}_{p(1)} \cdots u^{\ell+1}_{p(\ell+1)} = 1,
\end{gather*}
where the sum is over all permutations~$p$ of the set $\{1, \dots, \ell +1 \}$ and $\| p \|$ is the number of inversions
in~$p$.
The $*$-structure is def\/ined as in~\cite{projective}.

There is a~non-degenerate pairing $\langle \cdot, \cdot \rangle$ between ${\rm U}_{q}(\mathfrak{su}(\ell + 1))$ and
$\mathcal{A}({\rm SU}_{q}(\ell + 1))$, which is used to def\/ine the canonical left and right actions as $h \triangleright a=
a_{(1)} \langle h, a_{(2)} \rangle$ and $a \triangleleft h = \langle h, a_{(1)} \rangle a_{(2)}$, where we use
Sweedler's notation for the coproduct.
In this way $\mathcal{A}({\rm SU}_{q}(\ell + 1))$ can be seen as the algebra generated by the matrix coef\/f\/icients of the
f\/inite-dimensional representations of ${\rm U}_{q}(\mathfrak{su}(\ell + 1))$, similarly to the classical case.
The pairing can be extended in a~natural way to include also the action of the element $\hat{K}$ and its inverse.

There is a~faithful state on $\mathcal{A}({\rm SU}_{q}(\ell + 1))$, called the Haar state and denoted by~$\varphi$, which
generalizes the properties of the Haar integral in the classical case.
However, dif\/ferently from the classical case, the Haar state is not a~trace on $\mathcal{A}({\rm SU}_{q}(\ell + 1))$.
In particular its modular group is implemented by the element $K_{2\rho}$, in the sense that
\begin{gather}
\label{eq:haar-suq}
\varphi(ab) = \varphi(b K_{2\rho} \triangleright a\triangleleft K_{2\rho}).
\end{gather}

Consider now the left action of ${\rm U}_{q}(\mathfrak{su}(\ell + 1))$ on $\mathcal{A}({\rm SU}_{q}(\ell + 1))$ def\/ined~by
\begin{gather*}
\mathcal{L}_{h}a = a\triangleleft S^{-1}(h).
\end{gather*}
It can be used to def\/ine the coordinate algebra $\mathcal{A}(S_{q}^{2\ell+1})$ of the quantum sphere $S^{2\ell+1}_{q}$
as
\begin{gather*}
\mathcal{A}\big(S_{q}^{2\ell+1}\big) = \big\{a \in \mathcal{A}({\rm SU}_{q}(\ell + 1)): \mathcal{L}_{h}a = \varepsilon(h)a, \, \forall\, h \in
{\rm U}_{q}(\mathfrak{su}(\ell)) \big\}.
\end{gather*}
Finally, using the generator of ${\rm U}_{q}(\mathfrak{u}(1))$, which we denoted by $\hat{K}$, we def\/ine the coordinate
algebra $\mathcal{A}(\mathbb{C}{\rm P}_{q}^{\ell})$ of the quantum projective space $\mathbb{C}{\rm P}_{q}^{\ell}$ as
\begin{gather*}
\mathcal{A}\big(\mathbb{C}{\rm P}_{q}^{\ell}\big) = \big\{a \in \mathcal{A}\big(S_{q}^{2\ell+1}\big): \mathcal{L}_{\hat{K}}a = a\big\}.
\end{gather*}

Having def\/ined the coordinate algebra $\mathcal{A}(\mathbb{C}{\rm P}_{q}^{\ell})$, the next step in order to build a~spectral
triple is to introduce a~Hilbert space, on which elements of this algebra act as bounded ope\-ra\-tors.
Recall that the projective spaces $\mathbb{C}{\rm P}^{\ell}$ are only spin$^{c}$ manifolds when~$\ell$ is even.
Then one possibility is to complete the space of antiholomorphic forms, with the idea of def\/ining a~Dolbeault--Dirac
operator acting on it.
This is the strategy followed in~\cite{projectiveplane} for the case $\ell=2$, where a~$q$-analogue of the space of
antiholomorphic forms is introduced.

This strategy is generalized in~\cite{projective} for all quantum projective spaces.
We denote by $\Omega^{k}$ their~$q$-analogue of the space of antiholomorphic~$k$-form.
More generally, they also consider the possibility of twisting this module of~$k$-forms by a~line bundle $\Gamma_{N}$,
with the resulting space being denoted by $\Omega^{k}_{N}$, and with the space of forms corresponding to the case $N=0$.

The space $\bigoplus_{k=0}^{\ell} \Omega^{k}_{N}$
carries a~left action of ${\rm U}_{q}(\mathfrak{su}(\ell + 1))$ and can be
decomposed into irreducible representations.
The resulting decomposition takes the following form:
\begin{gather}
\Omega^{0}_{N}   \simeq \bigoplus_{m \in \mathbb{N}} V_{(m + c_{1}, 0, \dots, 0, m + c_{2})},
\nonumber
\\
\Omega^{k}_{N}   \simeq \bigoplus_{m \in \mathbb{N}} V_{(m + c_{3}, 0, \dots, 0, m + c_{4}) + e_{k}} \oplus V_{(m +
c_{5}, 0, \dots, 0, m + c_{6}) + e_{k+1}}
\qquad
\text{for}\quad 1 \leq k \leq \ell-1,
\nonumber
\\
\Omega^{\ell}_{N}   \simeq \bigoplus_{m \in \mathbb{N}} V_{(m + c_{7}, 0, \dots, 0, m + c_{8})}.
\label{forms}
\end{gather}
Here $c_{1}, \dots, c_{8}$ are integers depending on~$k$ and~$N$, but independent of~$m$.
The Hilbert spaces obtained as a~completions of $\bigoplus_{k=0}^{\ell} \Omega^{k}_{N}$ are denoted by $H_{N}$.

It is possible to introduce the~$q$-analogue of the Dolbeault operator, which we denote by $\bar{\partial}$, which maps
$\Omega^{k}_{N}$ into $\Omega^{k+1}_{N}$ and satisf\/ies $\bar{\partial}^{2} = 0$.
Similarly the adjoint $\bar{\partial}^{\dagger}$ maps $\Omega^{k+1}_{N}$ into $\Omega^{k}_{N}$ and satisf\/ies
$(\bar{\partial}^{\dagger})^{2} = 0$.
A~family of Dolbeault--Dirac operators, denoted by $D_{N}$ for $N \in \mathbb{Z}$, is def\/ined by taking suitable linear
combinations of $\bar{\partial}$ and $\bar{\partial}^{\dagger}$ on each $\Omega^{k}_{N}$.
The operator $D_{0}$ is the~$q$-analogue of the Dolbeault--Dirac operator on $\mathbb{C}{\rm P}^{\ell}$, while $D_{N}$ is the
twist of $D_{0}$ with the Grassmannian connection of a~certain line bundle.
In particular, if~$\ell$ is odd and $N = (\ell+1)/2$, then $D_{N}$ is the~$q$-analogue of the Dirac operator for the
Fubini--Study metric.

In this paper we do not need the precise form of $D_{N}$, but only an asymptotic form of its eigenvalues.
In particular, for our purposes, this turns out to be independent on the value of~$N$.
Using the decomposition~\eqref{forms}, it is possible to compute the eigenvalues of $|D_{N}|$ when restricted to the
space $\Omega^{k}_{N}$.
The information that we need is that these eigenvalues grow like $q^{-m}$ with $m \in \mathbb{N}$, see the discussion at
the end of~\cite{projective}.

\section{Modular spectral triples}
\label{sec:modular}
Consider now the restriction of the Haar state of $\mathcal{A}({\rm SU}_{q}(\ell + 1))$ to
$\mathcal{A}(\mathbb{C}{\rm P}_{q}^{\ell})$, denoted again by the symbol~$\varphi$.
It follows, using the def\/initions given in the previous section, that any element $a \in
\mathcal{A}(\mathbb{C}{\rm P}_{q}^{\ell})$ is invariant under the right action of $K_{2\rho}$, that is $a \triangleleft
K_{2\rho} = a$.
Therefore the modular property of the Haar state of $\mathcal{A}({\rm SU}_{q}(\ell + 1))$, given by~\eqref{eq:haar-suq},
becomes
\begin{gather*}
\varphi(ab) = \varphi(b K_{2\rho} \triangleright a),
\qquad
a,b \in \mathcal{A}\big(\mathbb{C}{\rm P}_{q}^{\ell}\big).
\end{gather*}
As we have remarked in the introduction, the non-commutative integral, def\/ined in the usual sense of spectral triples in
terms of $D_{N}$, does not coincide with the Haar state.
Indeed the former is a~trace while the latter is not.
This fact provides a~motivation to introduce a~twist in the def\/inition of the non-commutative integral, as we now
proceed to explain.

We denote by $K_{2\rho}$ the closure of the unbounded operator acting via the left action of $K_{2\rho}$ on
$\mathcal{A}({\rm SU}_{q}(\ell + 1))$, which is a~positive and invertible operator.
We assume for the moment that the operator $K_{2\rho} (D^{2}_{N}+1)^{-s/2}$ is trace-class for all $s > p$, for some $p
\geq 0$.
Let also assume that the following linear functional on $\mathcal{A}(\mathbb{C}{\rm P}_{q}^{\ell})$ makes sense
\begin{gather*}
\psi(a) = \limr{s}{p} \Tr\big(K_{2\rho} a\big(D_{N}^{2} + 1\big)^{-s/2}\big).
\end{gather*}
Then it can be shown, as will be done in Appendix~\ref{appendix}, that we have
\begin{gather}
\psi(ab)= \limr{s}{p} \Tr\big(K_{2\rho} ab \big(D_{N}^{2} + 1\big)^{-s/2}\big)
  = \limr{s}{p} \Tr\big(b K_{2\rho} a\big(D_{N}^{2} + 1\big)^{-s/2}\big)
\nonumber
\\
\phantom{\psi(ab)}
  = \limr{s}{p} \Tr\big(K_{2\rho} K_{2\rho}^{-1} b K_{2\rho} a\big(D_{N}^{2} + 1\big)^{-s/2}\big).
\label{eq:twist-formula}
\end{gather}
Here the non-trivial equality is the second one.
Then, since we have $K_{2\rho}^{-1} b K_{2\rho} = K_{2\rho}^{-1} \triangleright b$ for all $b \in
\mathcal{A}(\mathbb{C}{\rm P}_{q}^{\ell})$, we f\/ind the modular property
\begin{gather*}
\psi(ab) = \psi\big(K_{2\rho}^{-1} \triangleright b a\big),
\end{gather*}
which is equivalent to that of equation~\eqref{eq:haar-suq}.
Therefore in this way we should obtain a~linear functional on $\mathcal{A}(\mathbb{C}{\rm P}_{q}^{\ell})$ which has the
modular property of the Haar state.

Of course we should check that the assumptions made above are justif\/ied.
This will be done in the next sections, where we will check them explicitely in the case $a = 1$.
It turns out that this is enough to conclude in the general case, since~\cite[Lemma 2.1]{netu05} guarantees that if
$\psi(1)$ is well def\/ined then $\psi(a)$ is proportional to the Haar state.

Before getting into that, we should mention that this kind of construction f\/its into the framework of \textit{modular
spectral triples}, which was introduced in~\cite{modular1}.
The main motivation for this concept was to study algebras that do not admit non-trivial traces, but it was later
realized that it can also be used to study algebras which have faithful states but not faithful traces, as
in~\cite{modular2}.
The latter is the relevant case for this paper.
This notion has been formalized on the basis of examples where the modular group comes from a~circle action, so that
modif\/ications might be needed to handle more complicated examples.

\begin{definition}
Let $\mathcal{A}$ be a~unital $*$-subalgebra of~$N$, where~$N$ is a~semif\/inite von Neumann algebra acting on a~Hilbert
space~$H$.
Fix a~faithful normal strictly semif\/inite weight~$\phi$ with modular group $\sigma^{\phi}$.
We call the triple $(\mathcal{A},H,D)$ a~modular spectral triple if
\begin{enumerate}\itemsep=0pt
\item[1)] $\mathcal{A}$ is invariant under $\sigma^{\phi}$ and consist of analytic vectors for $\sigma^{\phi}$,
\item[2)] $D$ is a~self-adjoint operator af\/f\/iliated with the f\/ixed point algebra $N^{\sigma^{\phi}}$,
\item[3)] $[D,a]$ extends to a~bounded operator in~$N$ for all $a \in \mathcal{A}$,
\item[4)] $(D^2+1)^{-1/2}$ is compact with respect to the trace $\tau = \phi |_{N^{\sigma^{\phi}}}$.
\end{enumerate}
\end{definition}

We recall that a~semif\/inite weight~$\phi$ is \emph{strictly} semif\/inite if its restriction to $N^{\sigma^{\phi}}$ is
semif\/inite, and that the ideal of compact operators with respect to a~semif\/inite weight is the norm-closed ideal
generated by projections on which the weight takes f\/inite values.

It is worth noting that if $\mathcal{A}$ is pointwise invariant under the modular group $\sigma^{\phi}$ then we are
essentially back to the semif\/inite case.
This observation makes clear the fact that the f\/ixed point algebra plays an important role in this def\/inition.
However, in examples it might well be that no element of~$\mathcal{A}$ is invariant under the modular group.

Regarding summability, the notion of spectral dimension can be adapted straightforwardly to this case by replacing the
trace with the state or weight under consideration.

\begin{definition}
A~modular spectral triple $(\mathcal{A},H,D)$ is called f\/initely summable if there exists some $s_{0}>0$ such that $\phi
((D^{2}+1)^{-s_{0}/2}) < \infty$. In this case, we def\/ine the \textit{spectral dimension} as
\begin{gather*}
p = \inf \big\{s>0: \phi \big(\big(D^{2}+1\big)^{-s/2}\big) < \infty \big\}.
\end{gather*}
\end{definition}

A modif\/ication of this notion has appeared in~\cite{kaad}, by replacing the condition of boundedness of the commutator
with the analogue one for a~twisted commutator.
An interesting example that makes use of this condition is the one provided for ${\rm SU}_{q}(2)$ in~\cite{suq2-kaad}.

\section{Quantum dimension}
\label{sec:quantum}

Motivated by the previous section, we now want to introduce the tools needed to compute the spectral dimension of
$D_{N}$ with respect to the weight def\/ined by $\Tr(K_{2\rho} \cdot)$.
This computation is strictly related to the notion of quantum dimension, that we now review.

Given a~f\/inite-dimensional irreducible representation~$T$ of a~Drinfeld--Jimbo algebra ${\rm U}_{q}(\mathfrak{g})$, its
\textit{quantum dimension} is def\/ined as the number $\Tr(T(K_{2\rho}))$, where the trace is taken over the
vector space that carries the representation~$T$, see for example~\cite{KS}.
In the classical case, that is for $q=1$, the quantum dimension coincides with the dimension of the vector space.
In the context of quantum groups the notion of quantum dimension appears, for example, in the~$q$-analogue of the Schur
orthogonality relations.

In the classical case, if we consider a~f\/inite-dimensional representation of a~Lie algebra $\mathfrak{g}$ with highest
weight~$\Lambda$, the dimension of the associated vector space $V_{\Lambda}$ can be computed from the \textit{Weyl
dimension formula}, which reads as
\begin{gather*}
\dim V_{\Lambda} = \prod\limits_{\alpha>0}\frac{(\Lambda+\rho, \alpha)}{(\rho, \alpha)},
\end{gather*}
where the product is over the positive roots and~$\rho$ is the Weyl vector, def\/ined as the half-sum of the positive
roots.
There is also a~$q$-analogue of this formula, see~\cite{Fuc} and references within (notice that our~$q$ is $q^{1/2}$ in
the notation of the book).
It allows to compute the quantum dimension of a~representation with highest weight~$\Lambda$ as
\begin{gather*}
\dim_{q} V_{\Lambda} = \prod\limits_{\alpha>0} \frac{[(\Lambda + \rho, \tilde{\alpha})]}{[(\rho, \tilde{\alpha})]},
\end{gather*}
where we use the usual notion of~$q$-number
\begin{gather*}
[x] = \frac{q^{-x}-q^{x}}{q^{-1}-q}
\end{gather*}
and $\tilde{\alpha}=2\alpha/(\theta,\theta)$ where~$\theta$ is the highest root.
Note that an explicit normalization is needed for the positive roots, dif\/ferently from the classical case.

Our aim is now to compute the quantum dimension for any of the irreducible representations that appear in the
decomposition~\eqref{forms}.
More precisely we are only interested in the asymptotics of this value when $m \to \infty$, since this is the only
contribution that matters in the computation of the spectral dimension.

We need to review some facts about the root system of $\mathfrak{su}(\ell + 1)$, whose elements can be considered as
vectors in $\mathbb{R}^{\ell + 1}$.
The simple roots are given by $\alpha_{i} = e_{i} - e_{i+1}$ with $1 \leq i \leq \ell$.
The positive roots are given by $\alpha_{ij} = e_{i} - e_{j}$, with $1 \leq i < j \leq \ell+1$, and we note that they
can be written in terms of the simple roots as $\alpha_{ij} = \sum\limits_{k=i}^{j-1} \alpha_{k}$.
Their scalar product is $(\alpha_{ij}, \alpha_{ij}) = 2$.
In particular $(\theta, \theta) = 2$, so that $\tilde{\alpha} = \alpha$ in the Weyl formula.

We also need the basis of the fundamental weights, which we denote by $\omega_{i}$.
They are connected to the simple roots via the Cartan matrix~$A$ as $\alpha_{i} = \sum\limits_{j=1}^{\ell} A_{ij}
\omega_{j}$.
The fundamental weights are dual to the simple roots in the sense that
\begin{gather*}
\frac{2(\alpha_{i}, \omega_{j})}{(\alpha_{i}, \alpha_{i})} = \delta_{ij}.
\end{gather*}
Since in our case $(\alpha_{i}, \alpha_{i}) = 2$ this relation becomes $(\alpha_{i}, \omega_{j}) = \delta_{ij}$.
Finally we recall that the Weyl vector~$\rho$, which is usually def\/ined as the half-sum of the positive roots, can be
written in the basis of the fundamental weights in the simple form $\rho = \sum\limits_{j=1}^{\ell} \omega_{j}$.

In the following we will use the notation $f(m) \sim g(m)$ for $m \to \infty$ to mean that\linebreak
$\lim\limits_{m \to \infty} f(m) / g(m) = C$, where~$C$ is non-zero.
It is not dif\/f\/icult to determine~$C$ explicitly in the expressions we will present, but this is unnecessary for our
purposes.

\begin{proposition}
\label{prop:dimq}
Let $\Lambda = n_{1} \omega_{1} + n_{a} \omega_{a} + n_{\ell} \omega_{\ell}$ be a~dominant weight, where $n_{1} = m +
c_{1}$, $n_{\ell} = m + c_{2}$ with $m \in \mathbb{N}$, $c_{1}, c_{2} \in \mathbb{N}$ and $n_{a} = 0,1$ with $2 \leq a\leq \ell-1$.
Then for the corresponding quantum dimension we have $\dim_{q}(V_{\Lambda}) \sim q^{-2\ell m}$ for $m \to \infty$.
\end{proposition}
\begin{proof}
First observe that $[x] = (q^{x} - q^{-x})/(q - q^{-1}) \sim q^{-x}$ for $x \to \infty$, since we are assuming that $0 <
q < 1$.
Then we introduce the notation
\begin{gather*}
S_{i} = \prod\limits_{j=i+1}^{\ell+1}\frac{[(\Lambda + \rho,\alpha_{ij})]}{[(\rho, \alpha_{ij})]},
\end{gather*}
in such a~way that $\dim_{q}(V_{\Lambda})$ is given by the product of the $S_{i}$, that is
\begin{gather*}
\dim_{q}(V_{\Lambda}) = \prod\limits_{i=1}^{\ell}S_{i}.
\end{gather*}

Let us consider f\/irst the case $i=1$.
Using the formulae $(\alpha_{i}, \omega_{j}) = \delta_{ij}$ and $\alpha_{ij} = \sum\limits_{k=i}^{j-1} \alpha_{k}$ it is
immediate to show that we have
\begin{gather*}
(\Lambda,\alpha_{ij}) =
\begin{cases}
n_{1}, & j\leq a,
\\
n_{1}+n_{a}, & a<j<\ell+1,
\\
n_{1}+n_{a}+n_{\ell}, & j=\ell+1.
\end{cases}
\end{gather*}
Then for $m \to \infty$ we obtain
\begin{gather*}
S_{1} = \prod\limits_{j=2}^{\ell+1}\frac{[(\Lambda+\rho, \alpha_{ij})]}{[(\rho, \alpha_{ij})]} \sim q^{-(\ell - 1) m}
q^{- 2 m} = q^{-(\ell + 1) m}.
\end{gather*}
Similarly for $2\leq i\leq a$ we have
\begin{gather*}
(\Lambda,\alpha_{ij})=
\begin{cases}
0, & j\leq a,
\\
n_{a}, & a<j<\ell+1,
\\
n_{a}+n_{\ell}, & j=\ell+1.
\end{cases}
\end{gather*}
and for $m \to \infty$ we obtain
\begin{gather*}
S_{i} = \prod\limits_{j=i+1}^{\ell+1}\frac{[(\Lambda + \rho, \alpha_{ij})]}{[(\rho, \alpha_{ij})]} \sim q^{- m}.
\end{gather*}
Finally for $i\geq a+1$ we have
\begin{gather*}
(\Lambda,\alpha_{ij}) =
\begin{cases}
0, & j<\ell+1,
\\
n_{\ell}, & j=\ell+1,
\end{cases}
\end{gather*}
and for $m \to \infty$ we obtain
\begin{gather*}
S_{i} = \prod\limits_{j=i+1}^{\ell+1}\frac{[(\Lambda + \rho, \alpha_{ij})]}{[(\rho, \alpha_{ij})]} \sim q^{-m}.
\end{gather*}
Putting all together we f\/ind
\begin{gather*}
\dim_{q}(V_{\Lambda}) = S_{1} \left(\prod\limits_{i=2}^{a}S_{i}\right) \left(\prod\limits_{i=a+1}^{\ell}S_{i}\right)
 \sim q^{- (\ell + 1) m} \left(\prod\limits_{i = 2}^{a} q^{- m} \right) \left(\prod\limits_{i=a+1}^{\ell} q^{- m}\right)
\\
\phantom{\dim_{q}(V_{\Lambda})}
  = q^{- (\ell + 1) m} q^{- (a - 1) m} q^{- (\ell - a) m} = q^{- 2 \ell m}.\tag*{\qed}
\end{gather*}
\renewcommand{\qed}{}
\end{proof}

\section{Spectral dimension}
\label{sec:spectral}

Given the result of the previous section, it is now easy to prove the main result of the paper.

\begin{theorem}
\label{thm:zeta}
The operator $K_{2\rho} (D_{N}^{2}+1)^{-s/2}$ is trace-class for $s > 2\ell$, and the corresponding spectral dimension
$($in the sense of modular spectral triples$)$ is $2 \ell$.
Moreover
\begin{gather*}
\limr{s}{2\ell} \Tr\big(K_{2\rho} \big(D_{N}^{2}+1\big)^{-s/2}\big)
\end{gather*}
exists and is non-zero.
\end{theorem}
\begin{proof}

The Hilbert space $H_{N}$ is the completion of $\bigoplus_{k=0}^{\ell} \Omega^{k}_{N}$
and each $\Omega^{k}_{N}$ can be
decomposed into irreducible representations of ${\rm U}_{q}(\mathfrak{su}(\ell + 1))$ as in~\eqref{forms}.
As shown in~\cite{projective}, the operator $D_{N}^{2}$ restricted to the space $\Omega^{k}_{N}$ can be expressed in
terms of the Casimir operator of ${\rm U}_{q}(\mathfrak{su}(\ell + 1))$.
Therefore it acts as a~multiple of the identity in each irreducible representation.

The only representations which appear in the decomposition~\eqref{forms} are those of weight $(m + c_{1,k,N}, 0, \dots,
0, m + c_{2,k,N}) + e_{k}$, where $m \in \mathbb{N}$, $2 \leq k \leq \ell$ and $c_{1,k,N}$, $c_{2,k,N}$ are some positive
integers depending on~$k$ and~$N$.
We denote the vector space that carries such a~representation by $V_{m,k,N}$ and the corresponding eigenvalue of
$D_{N}^{2}$ by $\lambda_{m,k,N}^{2}$.
Finally denoting by $\Tr_{m,k,N}$ the trace on the vector space $V_{m,k,N}$ we have that
\begin{gather*}
\Tr_{m,k,N}\big(K_{2\rho}\big(D_{N}^{2} + 1\big)^{-s/2}\big) = \dim_{q}(V_{m,k,N}) \big(\lambda_{m,k,N}^{2} + 1\big)^{-s/2}.
\end{gather*}
From Proposition~\ref{prop:dimq} we know that $\dim_{q}(V_{m,k,N}) \sim q^{-2 \ell m}$ for $m \to \infty$.
Moreover we know from~\cite{projective} that $\lambda_{m,k,N} \sim q^{-m}$.
Then $\Tr_{m,k,N}(K_{2\rho}(D_{N}^{2}+1)^{-s/2}) \sim q^{(s- 2\ell)m}$.

Finally the trace can be written in the form
\begin{gather*}
\Tr\big(K_{2\rho}\big(D_{N}^{2} + 1\big)^{-s/2}\big)   = \sum\limits_{k=0}^{\ell} \sum\limits_{m=1}^{\infty}
\Tr_{m,k,N}\big(K_{2\rho}(D_{N}^{2} + 1)^{-s/2}\big)
  \sim \sum\limits_{m=1}^{\infty} q^{(s - 2\ell)m}.
\end{gather*}
The series $\sum\limits_{m=1}^{\infty} q^{(s - 2\ell)m}$ is absolutely convergent for $s > 2\ell$, from which it follows
that the spectral dimension is $2 \ell$.
Moreover we have $\sum\limits_{m=1}^{\infty} q^{(s - 2\ell)m} = q^{s - 2 \ell} / (1 - q^{s - 2 \ell})$, from which one
easily f\/inds that the limit exists and is non-zero.
\end{proof}

It is also possible to prove that $\Tr(K_{2\rho} (D_{N}^{2}+1)^{-s/2})$ extends to a~meromorphic function on the
complex plane, similarly to the classical case, but we will not show it here.

We now give a~few comments on this result.
As we mentioned in the introduction, in the classical case the computation of the spectral dimension hinges on the
balance between the growth of the eigenvalues of~$D$ and the growth of their multiplicities.
On the other hand, in the case of~$q$-deformations the eigenvalues of~$D$ grow like~$q$-numbers, therefore
exponentially, while their multiplicities only grow polynomially.
This has the consequence of giving a~spectral dimension equal to zero for the spectral triples
$(\mathcal{A}(\mathbb{C}{\rm P}_{q}^{\ell}), H_{N}, D_{N})$.
Roughly speaking, the ef\/fect of the weight $\Tr(K_{2\rho} \cdot)$ is to replace the multiplicities of the
eigenvalues with their~$q$-analogues, therefore restoring the balance in the computation.
Indeed it can be argued that in this context the notion of quantum dimension is the most natural one, as seen from its
role in the formulation of the quantum orthogonality relations.

The same result for the spectral dimension is obtained by considering $K_{2\rho}^{-1}$, as follows from a~general
property of the quantum dimension.

\begin{corollary}
\label{cor:twist}
The results of Theorem~{\rm \ref{thm:zeta}} remain valid if $K_{2\rho}$ is replaced by $K_{2\rho}^{-1}$.
\end{corollary}
\begin{proof}
This follows from the identity $\Tr(K_{2\rho}^{-1}) = \Tr(K_{2\rho})$, where the trace is taken on the
vector space of an irreducible f\/inite-dimensional representation, see~\cite[\S~7.1.6]{KS}.
We give an outline of the proof for the benef\/it of the reader.
There is an algebra automorphism~$\eta$ of ${\rm U}_{q}(\mathfrak{su}(\ell + 1))$ which is given on the generators as
\begin{gather*}
E_{i} \mapsto F_{\ell + 1 - i},
\qquad
F_{i} \mapsto E_{\ell + 1 - i},
\qquad
K_{i} \mapsto K_{\ell + 1 - i}^{-1},
\qquad
1 \leq i \leq \ell,
\end{gather*}
as can be checked directly from the def\/ining relations of ${\rm U}_{q}(\mathfrak{su}(\ell + 1))$.
Since $\eta(K_{2 \rho}) = K_{2 \rho}^{-1}$, the trace of $K_{2 \rho}$ in $V_{n}$ is equal to the trace of $K_{2
\rho}^{-1}$ in $V_{n}$ twisted by~$\eta$.
But the latter module is isomorphic to $V_{n}$ itself, since they are both simple f\/inite-dimensional modules with the
same highest weight.
This proves the claimed identity.
\end{proof}

This simple corollary is interesting in view of its possible applications to twisted Hochschild (co)homology,
see~\cite[\S~2.2]{SLq(2)} and references therein.
It is known that for quantum groups there is a~dimension drop in Hochschild homology: this means that, if~$G$ is
a~semisimple group and we denote by $\mathcal{A}(G_{q})$ the associated quantized algebra of functions, then we have
$H_{n}(\mathcal{A}(G_{q})) = 0$, where~$n$ denotes the classical dimension of~$G$.
On the other hand, by using twisted Hochschild homology, that is by twisting appropriately the notion of Hochschild
homology, it is possible to avoid this dimension drop.
This was observed f\/irst in~\cite{SLq(2)} for ${\rm SL}_{q}(2)$ by direct computation and then generalized in~\cite{dual-comp}
to the general case.

Similar results hold for quantum homogeneous spaces as the Podle\'{s} spheres, as shown by the computations
in~\cite{podles-hom}.
For results on a~more general class of quantum homogeneous spaces see~\cite{Kra12}.
For the standard Podle\'{s} sphere the dimension drop is avoided by considering the twist~$\vartheta_{P}^{-1}$, where~$\vartheta_{P}$ the modular group of the Haar state.
Then the volume form, being a~twisted cycle, will pair non-trivially with a~twisted cocycle with twist~$\vartheta_{P}^{-1}$.

In view of the results mentioned above, we expect that they continue to hold also for the projective spaces~$\mathcal{A}(\mathbb{C}{\rm P}_{q}^{\ell})$.
Therefore, if we denote by~$\vartheta$ the modular group in this case, we expect to avoid the dimension drop in homology
by twisting with~$\vartheta^{-1}$.
Then, in view of our results, we can write a~natural candidate for a~twisted cocycle that has a~chance of pairing
non-trivially with the volume form.

\begin{corollary}
The functional on $\mathcal{A}(\mathbb{C}{\rm P}_{q}^{\ell})^{\otimes(2\ell+1)}$ defined~by
\begin{gather*}
\tilde{\psi} (a_{0}, \dots, a_{2\ell}) = \limr{s}{2\ell} \Tr\big(K_{2\rho}^{-1} a_{0} [D_{N},a_{1}] \cdots
[D_{N},a_{2\ell}] \big(D_{N}^{2}+1\big)^{-s/2}\big)
\end{gather*}
is a~twisted cocycle with twist $\vartheta^{-1}$.
\end{corollary}
\begin{proof}
It follows from Corollary~\ref{cor:twist} that this functional is well-def\/ined.
That it is a~twisted cocycle with twist $\vartheta^{-1}$ follows from the twisted trace property shown in
equation~\eqref{eq:twist-formula}, with $K_{2\rho}$ replaced by $K_{2\rho}^{-1}$, and from standard computations.
\end{proof}

For the case of the Podle\'{s} sphere, it is shown in~\cite{res-form} that such a~twisted cocycle is indeed non-trivial,
when $D_{N}$ is taken to be the Dirac operator introduced in~\cite{DS-pod}.

\appendix

\section{The twisted trace property}
\label{appendix}
In this appendix we give a~proof, under suitable assumptions, of the equality appearing in
equation~\eqref{eq:twist-formula}.
The proof holds quite generally, so that we do not make any reference to quantum projective spaces.
We consider a~triple $(\mathcal{A}, H, D)$, where $\mathcal{A}$ is a~unital algebra acting as bounded operators on the
Hilbert space~$H$, that is $\mathcal{A} \subset B(H)$, and~$D$ is a~self-adjoint unbounded operator.
Moreover let $\Delta_{\phi}$ be a~positive invertible operator acting on~$H$.
The assumptions we make on this data are essentially those of a~modular spectral triple:
\begin{enumerate}\itemsep=0pt
\item[1)] $\Delta_{\phi} a\Delta_{\phi}^{-1} \in \mathcal{A}$ for any $a \in \mathcal{A}$,
\item[2)] $\Delta_{\phi}$ and $(D^{2}+1)^{-1/2}$ commute,
\item[3)] $[D,a]$ extends to a~bounded operator for any $a \in \mathcal{A}$.
\end{enumerate}
We also make the following summability assumptions:
\begin{enumerate}\itemsep=0pt
\item[1)] $\Delta_{\phi} \big(D^{2}+1\big)^{-s/2}$ is trace-class for all $s > p$, with $p \geq 0$,
\item[2)] $\limr{s}{p} \Tr\big(\Delta_{\phi} a\big(D^{2} + 1\big)^{-s/2}\big)$ exists for all $a \in \mathcal{A}$.
\end{enumerate}
Here by $\lim\limits_{s \to p^{+}}$ we mean the one-sided limit from the right.
We note in passing that these two summability conditions can be related to the semif\/inite theory, see~\cite{asym-zeta}.

Given these assumptions, we can def\/ine a~linear functional on $\mathcal{A}$~by
\begin{gather*}
\psi(a) = \limr{s}{p} \Tr\big(\Delta_{\phi} a\big(D^{2} + 1\big)^{-s/2}\big).
\end{gather*}

\begin{proposition}
With the same assumptions as above, the linear functional $\psi: \mathcal{A} \to \mathbb{C}$ satisfies the twisted trace
property $\psi(ab) = \psi\big(\Delta_{\phi}^{-1} b \Delta_{\phi} a\big)$ for all $a,b \in \mathcal{A}$.
\end{proposition}

\begin{proof}
The crucial step of the proof is to show that
\begin{gather*}
\psi(ab)= \limr{s}{p} \Tr\big(\Delta_{\phi} ab \big(D^{2} + 1\big)^{-s/2}\big)
  = \limr{s}{p} \Tr(\Delta_{\phi} a\big(D^{2} + 1\big)^{-s/2} b),
\end{gather*}
or equivalently that the following quantity vanishes
\begin{gather}
\label{eq:step2}
\limr{s}{p} \Tr\big(\Delta_{\phi} a\big[\big(D^{2} + 1\big)^{-s/2},b\big]\big) = 0.
\end{gather}
It is enough to consider $a = 1$ in this last equation, since using H\"older's inequality
\begin{gather*}
\big|\Tr\big(\Delta_{\phi} a\big[\big(D^{2} + 1\big)^{-s/2}, b\big]\big)\big| \leq \big\| \Delta_{\phi} a\Delta_{\phi}^{-1} \big\|
\Tr\big(\big|\Delta_{\phi} \big[\big(D^{2} + 1\big)^{-s/2}, b\big] \big|\big).
\end{gather*}

We proceed similarly to~\cite[Theorem 10.20]{green-book}, but taking care of the presence of the modular operator $\Delta_{\phi}$.
First of all we write $p = k \bar{r}$, with f\/ixed $k \in \mathbb{N}$ and $0 < \bar{r} < 1$ (notice that if~$p$ is an
integer we can set $k=2p$ and $\bar{r}=1/2$).
With this convention we can write any~$s$ in a~suf\/f\/iciently small neighbourhood of~$p$ as $s = kr$, for some $0<r<1$.

Then, using simple commutator identities, we obtain
\begin{gather*}
\big[\big(D^{2}+1\big)^{-s/2},b\big]= \sum\limits_{j=1}^{k} \big(D^{2}+1\big)^{-(j-1)r/2} \big[\big(D^{2}+1\big)^{-r/2}, b\big] \big(D^{2}+1\big)^{-(k-j)r/2}
\\
\phantom{\big[\big(D^{2}+1\big)^{-s/2},b\big]}
  = - \sum\limits_{j=1}^{k} \big(D^{2}+1\big)^{-jr/2} \big[\big(D^{2}+1\big)^{r/2}, b\big] \big(D^{2}+1\big)^{-(k-j+1)r/2}.
\end{gather*}
We introduce the notation
\begin{gather*}
R_{j} = \big(D^{2}+1\big)^{-jr/2} \big[\big(D^{2}+1\big)^{r/2}, b\big] \big(D^{2}+1\big)^{-(k-j+1)r/2}.
\end{gather*}
Let $p_{j}$ and $q_{j}$ be numbers such that $p_{j}^{-1} + q_{j}^{-1} = 1$.
Then we have
\begin{gather*}
\Delta_{\phi} R_{j}= \Delta_{\phi}^{p_{j}^{-1}} \Delta_{\phi}^{q_{j}^{-1}} \big(D^{2}+1\big)^{-jr/2}\Delta_{\phi}^{-q_{j}^{-1}} \Delta_{\phi}^{q_{j}^{-1}}
\\
\phantom{\Delta_{\phi} R_{j}=}
  \times \big[\big(D^{2}+1\big)^{r/2}, b\big] \Delta_{\phi}^{-q_{j}^{-1}} \Delta_{\phi}^{q_{j}^{-1}} \big(D^{2}+1\big)^{-(k-j+1)r/2}.
\end{gather*}
Since we assumed that~$D$ and $\Delta_{\phi}$ commute, this can be rewritten as
\begin{gather*}
\Delta_{\phi} R_{j} = \Delta_{\phi}^{p_{j}^{-1}} \big(D^{2}+1\big)^{-jr/2} \big[\big(D^{2}+1\big)^{r/2},
\Delta_{\phi}^{q_{j}^{-1}}b\Delta_{\phi}^{-q_{j}^{-1}}\big] \Delta_{\phi}^{q_{j}^{-1}} \big(D^{2}+1\big)^{-(k-j+1)r/2}.
\end{gather*}
Now from H\"older's inequality it follows that
\begin{gather*}
\Tr(|\Delta_{\phi} R_{j}|) \leq C_{j} \Tr\big(\Delta_{\phi} \big(D^{2}+1\big)^{-jp_{j}r/2}\big)^{p_{j}^{-1}}
\Tr\big(\Delta_{\phi} \big(D^{2}+1)^{-(k-j+1\big)q_{j}r/2}\big)^{q_{j}^{-1}},
\end{gather*}
where $C_{j} = \| [(D^{2}+1)^{r/2}, \Delta_{\phi}^{q_{j}^{-1}}b\Delta_{\phi}^{-q_{j}^{-1}}] \|$.
It follows from general arguments, which use the boundedness of $[D,a]$ for every $a \in \mathcal{A}$, that this
quantity is f\/inite, see~\cite[Lemma 10.17]{green-book}.

Now we want to choose $p_{j}$ and $q_{j}$ in such a~way that the operators $\Delta_{\phi} (D^{2}+1)^{-jp_{j}r/2}$ and
$\Delta_{\phi} (D^{2}+1)^{-(k-j+1)q_{j}r/2}$ are trace-class, which in turn would show that $\Delta_{\phi} R_{j}$ is
trace-class.
Since by assumption we have that $\Delta_{\phi} (D^{2} + 1)^{-s/2}$ is trace-class for all $s > p$, this implies the
inequalities $j p_{j} r > p$ and $(k - j + 1) q_{j} r > p$.
Let us set
\begin{gather*}
p_{j} = \frac{s}{r(j-1/2)},
\qquad
q_{j} = \frac{s}{r(k-j+1/2)},
\end{gather*}
and notice that they satisfy the equality $p_{j}^{-1} + q_{j}^{-1} = 1$, as they should.
Then it is immediate to see that the inequalities $jp_{j}r > p$ and $(k-j+1)q_{j}r > p$ are satisf\/ied for $s \geq p$
(notice the equality sign).
Therefore we have proven that $\Delta_{\phi} R_{j}$ is trace-class and, since
\begin{gather*}
\Delta_{\phi} \big[\big(D^{2}+1\big)^{-s/2}, b\big] = - \sum\limits_{j=1}^{k} \Delta_{\phi} R_{j},
\end{gather*}
the same is true for this operator.
This shows that the limit in~\eqref{eq:step2} is equal to zero and proves the equality claimed at the beginning.

The rest of the proof is now trivial.
Using the trace property we get
\begin{gather*}
\psi(ab)= \limr{s}{p} \Tr\big(\Delta_{\phi} a\big(D^{2} + 1\big)^{-s/2} b\big)
  = \limr{s}{p} \Tr\big(b \Delta_{\phi} a\big(D^{2} + 1\big)^{-s/2}\big)
\\
\phantom{\psi(ab)}
  = \limr{s}{p} \Tr\big(\Delta_{\phi} \Delta_{\phi}^{-1} b \Delta_{\phi} a\big(D^{2} + 1\big)^{-s/2}\big).
\end{gather*}
But this shows that $\psi(ab) = \psi(\Delta_{\phi}^{-1} b \Delta_{\phi} a)$, which concludes the proof.
\end{proof}

\subsection*{Acknowledgements}
I wish to thank Jens Kaad for helpful comments on a~f\/irst version of this paper.
I also want to thank the anonymous referees, whose observations have improved this presentation.

\pdfbookmark[1]{References}{ref}
\LastPageEnding

\end{document}